\documentclass[11 pt]{amsart}
\usepackage{graphicx, mathrsfs}
\usepackage{amssymb,amsmath,amsthm}
\usepackage{mathrsfs,dsfont}

\theoremstyle{plain}
\newtheorem{theorem}{Theorem}[section]

\newtheorem{proposition}{Proposition}[section]
\newtheorem{lemma}[proposition]{Lemma}
\newtheorem{defi}{Definition}[section]

\newtheorem*{coro}{Corollary}

\numberwithin{equation}{section} 

\newcommand{\RR}{\mathbb{R}}

\newcommand{\p}{\partial}

\newcommand{\CA}{\mathscr{A}}

\newcommand{\CD}{\mathcal{D}}
\newcommand{\CH}{\mathcal{H}}
\newcommand{\CM}{\mathcal{M}}
\newcommand{\TD}{\SD}

\newcommand{\D}{{\bf D}}
\newcommand{\R}{\mathcal{R}}
\newcommand{\CN}{\mathcal{N}}
\newcommand{\SD}{\mathscr{D}}
\newcommand{\SR}{\mathscr{R}}

\newcommand{\lf}{\left}
\newcommand{\rt}{\right}

\begin{document}

\title[Interface Problems]{A priori estimates for  fluid Interface problems}
\author[Shatah]{Jalal Shatah$^\dagger$}
\thanks{$^\dagger$ The first author is funded in part by NSF DMS
0203485.}
\address{$^\dagger$Courant Institute of Mathematical Sciences\\
251 Mercer Street\\ New York, NY 10012} \email{shatah@cims.nyu.edu}
\author[Zeng]{Chongchun Zeng$^*$}
\address{$^*$School of Mathematics\\
Georgia Institute of Technology\\ Atlanta, GA 30332}
\email{zengch@math.gatech.edu}
\thanks{$^*$ The second author is funded in part by NSF
DMS 0627842 and the Sloan Fellowship.}
\date{}
\begin{abstract}
We consider the regularity of an interface between two incompressible and inviscid fluids flows in the presence of surface tension.  We obtain local in time estimates on the interface in $H^{\frac32k +1}$ and the velocity fields  in $H^{\frac32k}$.   These estimates are obtained using geometric considerations which show that the Kelvin-Helmholtz instabilities are a consequence of a curvature calculation.
\end{abstract}
\maketitle

\section{Introduction}
In this manuscript we consider the interface problem between two incompressible and inviscid  fluids that occupy  domains $\Omega_t^+$ and $\Omega_t^-$ in $\RR^n$, $n \ge 2$, at time $t$.  We assume that $\RR^n = \Omega_t^+ \cup \Omega_t^-\cup S_t$ where $S_t = \p \Omega^\pm_t$, and let $p_\pm: \Omega_t^\pm \to \RR$, $v_\pm : \Omega^\pm_{t} \to \RR^n$, and
the constant $\rho_\pm>0$ denote the pressure, the velocity vector
field, and the density respectively.   On the interface $S_t$,  we let  $N_\pm(t, x)$, $x \in S_t$ denote the unit outward normal of $\Omega_t^\pm$ (thus  $N_+ + N_- =0$),   $H(t, x) \in (T_x S_t)^\perp$ denote the
mean curvature vector, and $\kappa_\pm = H \cdot N_\pm$. We also assume that there is surface tension  on the interface given by the mean curvature.  Thus  the free boundary problem for the Euler equation that we consider here  is given by
\begin{equation} \tag{E}\begin{cases}
\rho (v_t + \nabla_v v) = -\nabla p, \qquad & x\in \RR^n \smallsetminus
S_t \\
\nabla \cdot v =0, \qquad & x\in \RR^n \smallsetminus S_t,
\end{cases}\end{equation}
The boundary conditions for the interface evolution and the pressure are
\begin{equation}
\tag{BC}\begin{cases}
\p_t + v_\pm \cdot\nabla \quad \text{is tangent to} \quad
\bigcup\limits_{t}S_t\subset \mathbb{R}^{n+1},\\
p_+(t,x) - p_-(t,x) = \kappa_+ (t,x), \qquad x\in S_t.
\end{cases}
\end{equation}
where we introduced the notation
$v=v_+\mathds{1}_{\Omega_+} + v_-\mathds{1}_{\Omega_-}: \RR^n \smallsetminus S_t \to \RR^n$, etc.

The boundary conditions (BC) are a consequence of assuming that 1) the interface velocity is given by the normal component of the velocity $v_\pm^\perp=v_\pm\cdot N_\pm$,  and that 2) the surface tension on the interface $S_t$  is given by the mean curvature of the surface.  A weak formulation for the Euler flow with this form of surface tension is given by
\begin{equation*} \begin{cases}
(\rho v)_t + \nabla_v\rho  v = -\nabla p + H(t,x)\delta(S_t), \qquad &x\in \RR^n \\
\nabla \cdot v =0,& x\in \RR^n,
\end{cases}\\
\end{equation*}
where $\delta$ is the Dirac mass distribution.  This weak formulation implies the boundary condition for the pressure stated in (BC).

Here we consider the problem where $\Omega^+$ is compact and derive a priori estimates, local in time,  to the problem (E, BC) that prove bounds on $v(t, \cdot) \in H^{\frac32 k} (\RR^n \smallsetminus S_t)$ and $ S_t \in H^{\frac32 k+1}$ for ${\frac32 k} > \frac n2 +1$.  The assumption that $\Omega^+$ is compact is not necessary. In fact the same estimates hold if we assume that $S_t$ is either periodic or asymptotically flat.  A more interesting observation is that our proof works verbatim  for the case where there are several fluids occupying regions $\Omega_t^i$ with interfaces $S^i_t$. Thus we can treat more general setting than the existing literature.


The interface problem between two fluids has been studied extensively  in the math and physics literature.
In the absence of surface tension it is well known  that the interface problem between two inviscid and incompressible fluids is ill-posed due to the Kelvin-Helmholtz instability, and it is argued on physical basis that   the surface tension is a regularizing force that should make the problem well posed.  In~\cite{bhl93}, Beal, Hou, and Lowengrub demonstrated that the surface tension makes the linearized problem well-posed.   For the full nonlinear problem rigorous results  have been obtained  for irrotational velocities.  In this case   the problem can be completely reduced to the interface evolution with nonlocal operators. For this problem, Iguchi, Tanaka, and Tani~\cite{itt97} proved the local well-posedness in 2 dimensions with initial interface almost flat and initial velocity almost zero. For the general irrotational problem, Ambrose~\cite{am03} and, more recently, Ambrose and Masmoudi~\cite{am06} proved the local well-posedness in 2 and 3 dimensions, respectively.

We should note that for irrotational flows without surface tension the interface problem is given by he Birkhoff-Rott differential-integral equation.  Several results were obtained in this case such as those obtained by  Sulem,  Sulem,  Bardos,  and Frisch~\cite{ssbf81}, and Wu~\cite{wu06}.  We also note that without surface tension one can consider  the Euler equation with discontinuous velocity fields.  Although  the interface problem is ill-posed, due to  the Kelvin-Helmhotz instability,  weak solutions to the Euler equation which may include such discontinuities have been considered by DiPerna and Majda\cite{dm87}, Delort~\cite{de91}, and others.  There is also a rich literature of numerical studies of the interface problem, see for example~\cite{hls97} and references therein.

A related problem to the interface problem  is the water wave problem where there is only one fluid. In this case The Rayleigh-Taylor instability, instead of the Kelvin-Helmholtz instability, may occur. Such problems have been extensively studied and there is a vast literature on this subject, see for example, \cite{wu97}, \cite{wu99},
\cite{li05}.

Our approach in obtaining energy estimate for the interface problem  is similar to the water waves problem treated in \cite{sz06} in that it is geometric in nature.  It is based on the well known fact that these free boundary problems have a variational formulation on a subspace of volume  preserving homeomorphisms.  We use this variational approach to determine the terms that should be included in the energy. Of course these terms are identified as being the highest order terms of the linearized problem to (E, BC).  It is worth noting that from our analysis of the operators involved in the linearized
problem, the Kelvin-Helmholtz instability appears naturally as a
consequence of the negative semi-~definiteness of the leading part of
the unbounded curvature operator of the infinite dimensional
manifold of admissible Lagrangian coordinate maps. The surface
tension, created by the potential energy of the surface area,
generates a higher order positive operator that makes the linear problem well-posed and help to establish the
energy estimates. The well-posedness of the full problem well be addressed in a forthcoming article.

Our paper is organized as follows.  In section 2 we explain how to determine the pressure from the velocity.  In section 3 we give a variational formulation of the problem as a constrained variational problem for volume preserving maps. We use this formulation to motivate our definition of energy.  In section 4 we prove that our energy controls the Sobolev norm of the velocity and the mean curvature and derive bounds on the energy.  Some of the details in the geometric calculations are omitted since they are given in details in \cite{sz06} and are available as notes on the web at http://ww.math.gatech.edu/~zengch/notes/notes1.pdf.



\noindent {\bf Notation}   All notations will be defined as they are introduced.
In addition a list of symbols will be given at the end of the paper
for a quick reference.
The regularity of the domains $\Omega^\pm_t$ is characterized by the local regularity of $S_t$ as graphs. In general, an $m$-dimensional manifold $\CM \subset \RR^n$ is said to be of class
$C^k$ or $H^s$, $s > \frac m2$, if, locally in linear frames, $\CM$
can be represented by graphs of $C^k$ or $H^s$ mappings,
respectively.

 As   in~\cite{sz06}  $\Delta_\pm^{-1}$ denote the inverse Laplacian with zero Dirichlet data, $\CH_\pm$ denote the harmonic extension of functions defined on $S_t$ into $\Omega^\pm_t$, and
$\CN_\pm$ denote the Dirichlet  to Neuman operators in the domain $\Omega_t^\pm$. Given two fluids in $\Omega^\pm_t$ with constant densities $\rho_\pm$ we denote by
$
\CN =\frac 1{\rho_+}\CN_+ + \frac 1{\rho_-}\CN_-,
$
the operator $\CN^{-1}$ acts on function with mean zero and its range are also functions with mean zero.  For any quantity $q$  defined on  $\RR^n \smallsetminus S_t$  we write
$q=q_+\mathds{1}_{\Omega_+} +q_-\mathds{1}_{\Omega_-}$ where $q_\pm = q\mathds{1}_{\Omega^\pm}$.


\section{Determining the pressure}
In this section we explain how to express the pressure in terms of the velocity in this setting which is less clear than the  free boundary problem of water wave in  vacuum where the boundary conditions of $p_\pm$ are obvious.

To determine the boundary value of $p_\pm$ we take the dot product of Euler's equation (E)  with $N_\pm$
$$
-N_\pm \cdot \nabla p_\pm = \rho_\pm \D_{t\pm}(v_\pm\cdot N_\pm) -   \rho_\pm v_\pm\cdot  \D_{t\pm} N_\pm
$$
and using the fact that $v_+^\perp + v_-^\perp =0$, we obtain
$$
\frac 1{\rho_+} \nabla_{N_+} p_+ + \frac 1{\rho_-} \nabla_{N_-} p_- =v_+\cdot   \D_{t+} N_+ + v_-\cdot   \D_{t-} N_- - \nabla_{v_+^\top - v_-^\top} \,v_+^\perp.
$$

Substituting  the  formula for $\D_{t_\pm}
N_\pm$, which has been calculated in~\cite{sz06},
\begin{equation} \label{E:dtn}
\D_{t_\pm} N_\pm = - ((D v_\pm)^*(N_\pm))^\top
\end{equation}
we have
\[
\frac 1{\rho_+} \nabla_{N_+} p_+ + \frac 1{\rho_-} \nabla_{N_-} p_-
= \Pi_+(v_+^\top, v_+^\top) + \Pi_-(v_-^\top, v_-^\top) -
2\nabla_{v_+^\top - v_-^\top} \,v_+^\perp,
\]
where $\Pi_\pm$ is the second fundamental form of $S_t$ associated
to $N_\pm$, which satisfy $\Pi_+ +\Pi_-=0$. Since $p_\pm = \CH_\pm
(p_\pm|_{S_t}) + \Delta_\pm^{-1} \Delta p_\pm$ in $\Omega_t^\pm$, we
have
\[ \begin{split}
\frac 1{\rho_+} \CN_+ p_+ +\frac 1{\rho_-}\CN_- p_- = - \frac
1{\rho_+} \nabla_{N_+}& \Delta_+^{-1} \Delta p_+ - \frac 1{\rho_-}
\nabla_{N_-} \Delta_-^{-1} \Delta p_- \\
&+\Pi_+(v_+^\top, v_+^\top) + \Pi_-(v_-^\top, v_-^\top)
-2\nabla_{v_+^\top - v_-^\top} \,v_+^\perp \qquad \text{ on } S_t.
\end{split}\]
The boundary condition $p_+ - p_- = \kappa_+$ on $S_t$ stated in
(BC) implies that on $S_t$
\[\begin{split}
p_\pm = \CN^{-1} ( - \frac 1{\rho_\mp}\CN_\mp \kappa_\mp - \frac
1{\rho_+} \nabla_{N_+} \Delta_+^{-1} \Delta p_+ &- \frac 1{\rho_-}
\nabla_{N_-} \Delta_-^{-1} \Delta p_- )\\
& + \Pi_+(v_+^\top, v_+^\top) + \Pi_-(v_-^\top, v_-^\top)
-2\nabla_{v_+^\top - v_-^\top} \,v_+^\perp).
\end{split}\]
Finally, since $\nabla\cdot v =0$ in $\RR^n \smallsetminus S_t$, we have
from (E)
\begin{equation} \label{E:pressure1}
- \Delta p = \rho \nabla \cdot (\nabla_v v) = \rho \text{tr}(Dv)^2,
\qquad x \in \RR^n \smallsetminus S_t.
\end{equation}
Therefore,
\begin{equation} \label{E:pressure2} \begin{split}
p_\pm|_{S_t} = \CN^{-1} ( - \frac 1{\rho_\mp}\CN_\mp \kappa_\mp +&
\nabla_{N_+} \Delta_+^{-1} \text{tr}(Dv)^2 + \nabla_{N_-}
\Delta_-^{-1} \text{tr}(Dv)^2 ),\\
&+ \Pi_+(v_+^\top, v_+^\top) + \Pi_-(v_-^\top, v_-^\top)
-2\nabla_{v_+^\top - v_-^\top} \,v_+^\perp).
\end{split}\end{equation}
One can verify that the quantity $\CN^{-1}$ acts on in the above
has zero mean on $S_t$ and thus $p$ is well defined
by~\eqref{E:pressure1} and ~\eqref{E:pressure2}.

\section{Lagrangian formulation and the energy}
\label{S:geofunc}
This section is intended to explain the intuition behind the energy. It illustrates how to isolate the leading order  nonlinear terms $\CA$ and $\SR_0$ defined in  \eqref{E:CA1} and \eqref{E:R0} respectively.

In his 1966 seminal paper \cite{ar66}, V. Arnold pointed out that the Euler equation for an incompressible inviscid fluid can be viewed as the geodesic equation on the group of volume preserving diffeomorphisms. This point of view has been adopted and developed by several authors  such as D. G. Ebin and G. Marsden~\cite{em70}, A. Shnirelman~\cite{sh85}, and Y. Brenier \cite{br99}, to mention a few,  in their work on Euler's equations on fixed domains.  It is this point of view that  we adopted  to explain the motivation for our definition of energy for the water waves  problem \cite{sz06}, and it is this same  point of view that forms our starting point  to determine the appropriate energy for the interface problem .

\subsection{Lagrangian formulation of the problem} \label{SS:Lag}

Conservation of energy can be obtained from multiplying the Euler's equation (E)
by $v$, integrating on $\RR^n \smallsetminus S_t$, and using (BC) to
obtain the conserved energy $E_0$:
\begin{equation} \label{E:Econserv1}
E_0 = E_0(S_t, v)= \int_{\RR^n} \frac {\rho |v|^2}2 dx + \int_{S_t}
dS \triangleq \int_{\Omega_t} \frac {\rho |v|^2}2 dx + S(S_t),
\end{equation}
where $S(\cdot)$ denotes the surface area.

Let $u_\pm(t, y)$, $y \in \Omega_0^\pm$, be the Lagrangian
coordinate map solving
\begin{equation} \label{E:Lag1}
\frac {dx}{dt} = v(t, x), \qquad x(0)=y,
\end{equation}
then we have $v = u_t \circ u^{-1}$, and for any vector field $w$ on $x\in \RR^n \smallsetminus S_t$,
$\D_t w  = (w\circ u)_t \circ u^{-1}$.
Therefore in Lagrangian coordinates  the Euler's equation takes the form
\begin{equation} \label{E:euler1}
\rho u_{tt} = - (\nabla p) \circ u, \qquad u(0)= id_{\Omega_0},
\end{equation}
where the pressure $p$ is given by~\eqref{E:pressure1}
and~\eqref{E:pressure2}.

Since $v(t, \cdot)$ is divergence free in $\RR^n \smallsetminus S_t$,
then $u_\pm(t, \cdot)$ are volume preserving. Moreover, while $u_+(t,
\cdot)|_{S_0} = u_- (t, \cdot)|_{S_0}$ may not hold, it is clear
that $u_+ (t, S_0) = u_- (t, S_0)$.   Thus the Lagrangian coordinates  maps satisfy:\\
1) $\Phi_\pm:\bar{\Omega}^\pm \to  \Phi_\pm(\bar{\Omega}^\pm)$  a volume preserving homeomorphism.\\
2) $S\triangleq  \p \Phi_\pm(\Omega^\pm ) = \Phi(\p \Omega^\pm)$\\
Define
$$\Gamma = \{ \Phi = \Phi_+\mathds{1}_{\Omega^+} +  \Phi_-\mathds{1}_{\Omega^-};  \quad \Phi_\pm \text{ satisfy 1 and 2 above}\}.
$$
As a manifold, the tangent space of $\Gamma$ is given by divergence
free vector fields with matching normal component in  Eulerian
coordinates:
\[
T_\Phi \Gamma = \{ \bar w: \RR^n\smallsetminus S_0 \to \RR^n \mid
\nabla\cdot w=0 \text{ and } w_+^\perp + w_-^\perp|_{\Phi(S_0)} =0,
\text{where } w= (\bar w \circ \Phi^{-1}) \}.
\]
Here as in \cite{sz06} we are following the
convention that  for any vector field $X \, : \Phi(\Omega_0)\to\RR^n$ its
description in Lagrangian coordinates is given by $\bar X =
X\circ\Phi$.

Writing  $S(\Phi) =\int_{\Phi(S_0)} dS$ for the surface area of $\Phi(S_0)$, the
energy $E_0$ in Lagrangian coordinates can be written as:
\begin{equation} \label{E:Econserv2}
E_0 = E_0 (u, u_t) = \frac 12 \int_{\RR^n \smallsetminus S_0} \rho
|u_t|^2 dy + S(u), \qquad (u, u_t) \in T \Gamma
\end{equation}
where the volume preserving property of $u$ is used. This
conservation of energy suggests: 1) $T\Gamma$ be endowed with the
$L^2(\rho dy)$ metric\footnote{Including the density in the volume element $\rho dy$ introduces a factor of $\frac 1 \rho$ in front of the physical pressure.}; and 2) the free boundary problem of the
Euler's equation has a Lagrangian  action
\[
I(u) = \int \int_{\RR^n \smallsetminus S_0} \frac {\rho |u_t|^2}2 dy dt
- \int S(u) dt, \qquad u(t, \cdot) \in \Gamma.
\]
Let $\bar \SD$ denote the covariant derivative associated with the
metric on $\Gamma$, then a   critical path $u(t, \cdot)$ of $I$
satisfies
\begin{equation} \label{E:critical}
\bar \SD_t u_t + S' (u) =0.
\end{equation}

In order to verify that the Lagrangian coordinate map $u(t, \cdot)$
satisfying  (E) and (BC) is indeed a critical path of $I$, it
is convenient to calculate $\bar \SD$ and $S'$ by viewing $\Gamma$
as a submanifold of the Hilbert space $L^2(\RR^n \smallsetminus S_0,
\rho dy, \RR^n)$.

\noindent {\bf $(T_\Phi \Gamma)^\perp$ and  orthogonal decomposition
of vector fields}. For any vector field $X$ defined on $\Phi(\RR^n
\smallsetminus S_0)$,  Hodge decomposition suggests that we
decompose   $X$ into $X= w - \nabla \psi$,
with $\psi = \psi_+\mathds{1}_{\Omega_+} + \psi_-\mathds{1}_{\Omega_-}
$, so that $\bar w = w \circ \Phi \in
T_\Phi \Gamma$ and $\nabla \psi \circ \Phi \in (T_\Phi
\Gamma)^\perp$. For any $\bar Y \in T_\Phi \Gamma$, the
orthogonality $\int \nabla \psi  \cdot Y \rho dy =0$ implies
\[
\int_{\Phi(S_0)} \rho_+ \psi_+ Y_+^\perp + \rho_- \psi_- Y_-^\perp
dS =0.
\]
Therefore, $\psi$ must satisfy $\rho_+ \psi_+ = \rho_- \psi_-
\triangleq \psi^S$ on $\Phi(S_0)$.  This  suggests that, for any $\Phi
\in \Gamma$,
\[
(T_\Phi \Gamma)^\perp = \{-(\nabla \psi) \circ \Phi \mid \rho_+
\psi_+ = \rho_- \psi_- \; \text{ on } \; \Phi(S_0) \}.
\]
To prove this claim, we only need to find such a $\psi$ given $X$.
 From $X
= w -\nabla \psi$ and $w_+^\perp + w_-^\perp=0$, we have
\begin{equation} \label{E:BC2}
X_+^\perp + X_-^\perp = - \nabla_{N_+} \psi_+ - \nabla_{N_-} \psi_-
= - \CN \psi^S - \nabla_{N_+} \Delta_+^{-1} \Delta \psi -
\nabla_{N_-} \Delta_-^{-1} \Delta \psi.
\end{equation}
Since $\nabla \cdot w=0$, we obtain
\begin{equation} \label{E:normal1} \begin{cases}
-\Delta \psi = \nabla\cdot X \\
\psi_\pm|_{\Phi(S_0)} = \frac 1{\rho_\pm} \psi^S = - \frac
1{\rho_\pm} \CN^{-1} \lf( X_+^\perp + X_-^\perp - \nabla_{N_+}
\Delta_+^{-1} \nabla\cdot X - \nabla_{N_-} \Delta_-^{-1} \nabla\cdot
X\rt).
\end{cases}\end{equation}
It is easy to verify that $w = X - \nabla \psi$ satisfies $\bar
w = w \circ \Phi \in T_\Phi \Gamma$.

\noindent {\bf Computing $\bar \SD_t$ and $II_\Phi$}. Given a path
$u(t, \cdot) \in \Gamma$ and $\bar v= u_t$. Let $S_t = u(t,S_0)$.
Suppose $\bar w(t, \cdot) \in T_{u(t)} \Gamma$, then the covariant
derivative $\bar \SD_t \bar w$ and the second fundamental form
$II_{u(t)} (\bar w, \bar v)$ satisfy
\[
\bar w_t = \bar \SD_t \bar w + II_{u(t)} (\bar w, \bar v), \qquad
\bar \SD_t \bar w \in T_{u(t)} \Gamma, \quad II_{u(t)} (\bar w, \bar
v) \in (T_{u(t)} \Gamma)^\perp.
\]
Let $v = u_t \circ u^{-1} = \bar v \circ u^{-1}$ and $w=\bar w \circ
u^{-1}$ be the Eulerian coordinates description of $u_t$ and $w$, then for $X=\D_t w$
there exists $p_{v,w} =p_ {w,v}^+\mathds{1}_{\Omega^+} + p_{w,v}^-\mathds{1}_{\Omega^-} : \RR^n
\smallsetminus u(t, S_0): \to \RR$
determined by~\eqref{E:normal1}  such that
\begin{equation} \label{E:barCD}
\rho_+ p_{v,w}^+ = \rho_- p_{v,w}^- \text{ on } S_t, \quad II (\bar
w, \bar v) = - (\nabla p_{w,v})\circ u \in (T_{u(t)} \Gamma)^\perp.
\end{equation}
The Eulerian coordinates description of the covariant derivative is given by
\begin{equation} \label{E:barCD1}
\TD_t w = (\bar \SD_t \bar w) \circ u^{-1} = \D_t w + \nabla
p_{w,v}.
\end{equation}
The terms involving $X=\D_tw$ in equation \eqref{E:normal1} are expressed  as follows.
From $w_+^\perp + w_-^\perp=0$ on $u(t, S_0)$ and
identity~\eqref{E:dtn} we have
\begin{equation} \label{E:BC3}\begin{split}
(\D_{t_+} w_+) \cdot N_+ + (\D_{t_-} w_-)& \cdot N_- = \D_{t_+}
w_+^\perp + \D_{t_-} w_-^\perp + \nabla_{w_+^\top} v_+ \cdot N_+ +
\nabla_{w_-^\top} v_- \cdot N_-\\
=&\nabla_{v_+^\top - v_-^\top} w_+^\perp + \nabla_{w_+^\top -
w_-^\top} v_+^\perp - \Pi_+(v_+^\top, w_+^\top) - \Pi_-(v_-^\top,
w_-^\top).
\end{split}\end{equation}
And since
$
\nabla \cdot \D_t w = \text{tr}(DvDw)
$
then  $p_{w,v}$ is given by
\begin{equation} \label{E:II} \begin{cases}
-\Delta p_{w,v} = \text{tr}(DvDw)\\
p_{w,v}^\pm|_{S_t} = \frac 1{\rho_\pm} p_{w,v}^S = - \frac
1{\rho_\pm} \CN^{-1} \{ \nabla_{v_+^\top - v_-^\top} w_+^\perp +
\nabla_{w_+^\top - w_-^\top} v_+^\perp - \Pi_+(v_+^\top, w_+^\top) \\
\qquad \qquad \qquad\quad\; - \Pi_-(v_-^\top, w_-^\top) -
\nabla_{N_+} \Delta_+^{-1} \text{tr}(DvDw) - \nabla_{N_-}
\Delta_-^{-1} \text{tr}(DvDw)\}
\end{cases}\end{equation}
A more useful way to express the boundary value $p^S_{w,v}$ is as follows.
From the divergence decomposition formula
\begin{equation} \label{E:divB}
0 = \nabla\cdot v_\pm = \CD\cdot v_\pm^\top + \kappa_\pm v_\pm^\perp
+ N_\pm \cdot \nabla_{N_\pm} v_\pm \qquad \text{ on } S_t,
\end{equation}
where $\CD$ is the covariant derivative on $S_t$,  implies
\begin{equation} \label{E:XYN}
\nabla_{w_\pm} v_\pm \cdot N_\pm = \nabla_{w_\pm^\top} v_\pm \cdot
N_\pm  - \kappa_\pm w_\pm^\perp v_\pm^\perp - \CD \cdot (w_\pm^\perp
v_\pm^\top) + \nabla_{v_\pm^\top} w_\pm^\perp.
\end{equation}
Thus, we have
\begin{equation} \label{E:pS}\begin{split}
p_{w,v}^S = -\CN^{-1} \{ \nabla_{w_+} v_+ \cdot N_+ &+ \nabla_{w_-}
v_- \cdot N_- + \CD \cdot (w_+^\perp (v_+^\top - v_-^\top)) \\
& - \nabla_{N_+} \Delta_+^{-1} \text{tr}(DvDw) - \nabla_{N_-}
\Delta_-^{-1} \text{tr}(DvDw) \}.
\end{split} \end{equation}
Moreover, for any smooth function $f$ defined on $S_t$,  we have from the
Divergence theorem,
\[
\int_{S_t} - f \nabla_{N_\pm} \Delta_\pm^{-1} \text{tr} (Dv_\pm
Dw_\pm) \; dS = - \int_{\Omega_t^\pm} \nabla f_{\CH_\pm} \cdot
\nabla \Delta_\pm^{-1} \text{tr} (Dv_\pm Dw_\pm) + f_{\CH_\pm}
\text{tr} (Dv_\pm Dw_\pm) dx
\]
Again, by the Divergence Theorem  the first term integrates to zero and  the second term can be written as
\begin{equation} \label{E:pww1}\begin{split}
\int_{S_t} - f \nabla_{N_\pm} \Delta_\pm^{-1} \text{tr} (Dv_\pm
Dw_\pm) dS = &\int_{S_t} - f \nabla_{w_\pm} v_\pm \cdot N_\pm +
w_\pm^\perp \nabla f_{\CH_\pm} \cdot v_\pm dS \\
&- \int_{\Omega_t^\pm} D^2 f_{\CH_\pm} (v_\pm, w_\pm) dx.
\end{split}\end{equation}
Thus, using the decomposition $\nabla f_{\CH_\pm} = \nabla^\top f +
(\CN_\pm f) N_\pm$ and letting $f = -\CN^{-1} g$, we obtain
\begin{equation} \label{E:pww2}
\int_{S_t} g p_{w,v}^S dS = \int_{S_t} - w_+^\perp v_+^\perp (\CN_+
+ \CN_-) \CN^{-1}g\, dS + \int_{\RR^n \smallsetminus S_t} D^2 (\CH_\pm
(\CN^{-1} g)) (v, w) dx.
\end{equation}

\noindent {\bf Computing $S'(\Phi)$}. By the variation of surface
area formula, for any $\bar w \in T_\Phi \Gamma$ we have
\[
< S'(\Phi), \bar w>_{L^2(\RR^n \smallsetminus S_0, \rho dy)}
=\int_{\Phi(S_0)} \kappa_+ w_+^\perp \;dS = \int_{\Phi(S_0)}
\kappa_- w_-^\perp dS.
\]
We need to find the unique representation $S'(\Phi)$ in $T_\Phi
\Gamma$ of the above functional.
\begin{lemma} \label{L:SPrime}
For any smooth function $f_0: \Phi(S_0) \to \RR$, let $f_\pm = \pm
\frac 1{\rho_+ \rho_-} \CH_\pm \CN^{-1} \CN_\mp f_0$ and
$f=f_+ \mathds{1}_{\Omega^+} + f_-\mathds{1}_{\Omega^-}$, then we have $\nabla f \in T_\Phi \Gamma$ and for any $\bar w
\in T_\Phi \Gamma$,
\[
\int_{\Phi(S_0)} f_0 w_+^\perp dS = \int_{\RR^n \smallsetminus
\Phi(S_0)} w \cdot \nabla f\rho dx.
\]
\end{lemma}
The verification of the lemma is straightforward. Therefore, we have
\begin{equation}\label{E:SPrime}
S' (\Phi) = \nabla p_\kappa \qquad \text{ where } p_\kappa^\pm =
\frac 1{\rho_+ \rho_-} \CH_\pm \CN^{-1} \CN_\mp \kappa_\pm.
\end{equation}

\noindent {\bf Splitting of the pressure.} From~\eqref{E:pressure2},
\eqref{E:normal1} and~\eqref{E:SPrime}, it is clear that $\rho
(p_{v,v} + p_\kappa) = p$. Therefore, we obtain the well known
equivalence between equation~\eqref{E:critical} for critical paths
of $I$ and the Euler's equation (E) with the free boundary condition
(BC). The Euler' equation can also be written as
\begin{equation}\label{E:euler}
\D_t v + \nabla p_{v,v} + \nabla p_\kappa =0.
\end{equation}
Notice that the pressure splits into two terms, the first $p_{v,v}$
is the Lagrange multiplier, and the second $p_\kappa$ is due to
surface tension. These two terms will be treated differently in the
energy estimates.

\subsection{Linearization} \label{SS:linearization}

In order to analyze the free boundary problems of the Euler's
equation, it is natural to start with the linearization. The
Lagrangian formulation provides a convenient frame work for this
purpose. From~\eqref{E:critical}, the linearized equation is
\begin{equation} \label{E:criticalL}
\bar \SD_t^2 \bar w + \bar \SR (u_t, \bar w) u_t + \bar \SD^2 S(u)
(\bar w) =0, \qquad \bar w(t, \cdot) \in T_{u(t, \cdot)} \Gamma,
\end{equation}
where $\bar \SR$ is the curvature tensor of the infinite dimensional
manifold $\Gamma$.  Below we calculate $\bar \SR$ and $\bar \SD^2 S(u)$, which is  a linear operator on $T_u \Gamma$. Since these operators are self-adjoint, we will compute their  quadratic forms.

\noindent {\bf Computing $\bar \SD^2 S(u)$}.  The formula for $\bar \SD^2 S(u)$  was given in~\cite{sz06}

\begin{equation*}
\begin{split}
 \bar \SD^2 S(u) (\bar w, \bar w) = &\int_{S_t} \kappa_\pm
w_\pm^\perp (\kappa_\pm w_\pm^\perp + \CD \cdot w_\pm^\top) -
\kappa_\pm \nabla_{N_\pm} p_{w,w}^\pm - \kappa_\pm
\nabla_{w_\pm^\top} w_\pm \cdot N_\pm \\
&+ w_\pm^\perp\lf(-\Delta_{S_t}
w_\pm^\perp - w_\pm^\perp |\Pi|^2 + \nabla_{w_\pm^\top}
\kappa_\pm\rt) dS,
\end{split}
\end{equation*}
%
for any $\bar w \in T_u \Gamma$, where $\CD$ is the Riemannian
connection and $\Delta_{S_t}$ is the Beltrami-Lapalacian operator on
$S_t$. Of course $\bar \SD^2 S(u) (\bar w,
\bar w)$ is independent of the choice of the $+$ or $-$ sign.

Needless to say that this is a very complicated expression for $\bar
\SD^2 S(u) (\bar w, \bar w)$. We will single out its leading order
part. Since the value of $\bar{\SD}^2 S(u)$ does not depend on the choice
of $+$ or $-$ sign, we compute with the $+$
sign and assume that $S_t$ is a sufficiently smooth hypersurface. From
the Divergence Theorem,
\[\begin{split}
&|\bar \SD^2 S(u) (\bar w, \bar w) - \int_{S_t}\ |\nabla^\top
w_+^\perp|^2 \, dS|\\
\le & |\int_{S_t} \kappa_+ (\nabla_{N_+} p_{w,w}^+ +
\nabla_{w_+^\top} w_+ \cdot N_+ + \nabla_{w_+^\top} w_+^\perp) \;
dS| + C |w|_{L^2 (S_t)}^2.
\end{split}\]
To estimate the integral on the right side, we use the splitting of
$\nabla_{N_+} p_{w,w}^+$ on $S_t$
\[
\nabla_{N_+} p_{w,w}^+ = \CN_+(p_{w,w}^+|_{S_t}) + \nabla_{N_+}
\Delta_+^{-1} \Delta p_{w,w}^+ = \frac 1{\rho_+} \CN_+ p_{w,w}^S -
\nabla_{N_+} \Delta_+^{-1} \text{tr} (Dw)^2.
\]
and identities~\eqref{E:pww1} and~\eqref{E:XYN} to obtain
\[
|\bar \SD^2 S(u) (\bar w, \bar w) - \int_{S_t}\ |\nabla^\top
w_+^\perp|^2 \, dS| \le  |\int_{S_t} p_{w,w}^S \frac 1{\rho_+} \CN_+
\kappa_+  dS| + C (|w|_{L^2 (S_t)}^2 + |w|_{L^2 (\RR^n \smallsetminus
S_t)}^2).
\]
Finally, from~\eqref{E:pww2}, we have
\begin{equation} \label{E:D^2S}
|\bar \SD^2 S(u) (\bar w, \bar w) - \int_{S_t}\ |\nabla^\top
w_+^\perp|^2 \; dS| \le C (|w|_{L^2 (S_t)}^2 + |w|_{L^2 (\RR^n
\smallsetminus S_t)}^2).
\end{equation}
Using Lemma~\ref{L:SPrime}, we can define a self-adjoint
positive semi-~definite operator $\bar \CA(u)$ on $T_u \Gamma$ as the
leading order part of $\bar \SD^2 S(u)$. In the Eulerian
coordinates, $\bar \CA(u)$ takes the form
\begin{equation} \label{E:CA}
\CA (u) (w) =\nabla f_+ \mathds{1}_{\Omega^+} + \nabla f_-\mathds{1}_{\Omega^-}, \qquad \text{
where } f_\pm = \frac 1{\rho_+ \rho_-} \CH_\pm \CN^{-1} \CN_\mp
(-\Delta_{S_t}) w_\pm^\perp.
\end{equation}
Clearly $\bar \CA(u)$ is self-adjoint and satisfies
\begin{equation} \label{E:CA1}
\bar \CA(\bar w, \bar w) = \int_{S_t}\ |\nabla^\top w_+^\perp|^2 \,
dS.
\end{equation}
which is like a third order differential operator on  $\RR^n\smallsetminus S_t$. From~\eqref{E:D2S}, we can write
\begin{equation} \label{E:D2S}
\bar \SD^2 S(u) = \bar \CA(u) + \text{ at most 1st order diff.
operators}
\end{equation}

\noindent {\bf Computing $\bar \SR$}. 
For any $\bar v, \bar w \in
T_u \Gamma$, let $v = \bar v \circ u^{-1}$ and $w = \bar w \circ
u^{-1}$ we have
\[
\bar \SR(u) (\bar v, \bar w)\bar v \cdot \bar w = II_u (\bar v, \bar
v) \cdot II_u (\bar w, \bar w) - II_u(\bar v, \bar w)^2 =
\int_{\RR^n \smallsetminus S_t} \rho \nabla p_{v,v} \nabla p_{w,w} -
\rho |\nabla p_{v,w}|^2 \; dx.
\]
Assuming  that the hypersurface $S_t$ and $\bar v$ are
sufficiently smooth we single out the leading order term of  $\bar \SR(u) (\bar v, \cdot)\bar v$.

We first estimate $II_u (\bar v, \bar v) \cdot II_u (\bar w, \bar
w)$. From the Divergence Theorem and~\eqref{E:II},
\[
\int_{\RR^n \smallsetminus S_t} \rho \nabla p_{v,v} \nabla p_{w,w} dx =
\int_{S_t} p_{v,v}^S (\nabla_{N_+} p_{w,w}^+ + \nabla_{N_-}
p_{w,w}^-)\, dS + \int_{\RR^n\smallsetminus S_t} \rho p_{v,v} \text{tr}
(Dw)^2 \, dx
\]
Using~\eqref{E:BC2} and~\eqref{E:BC3} for the above boundary integral and applying the Divergence
Theorem twice for the interior integral, we obtain
\[\begin{split}
\int_{\RR^n \smallsetminus S_t}& \rho \nabla p_{v,v} \nabla p_{w,w} dx =
\int_{\RR^n \smallsetminus S_t} \rho D^2 p_{v,v} (w, w) \, dx +
\int_{S_t} p_{v,v}^S \{ -2\nabla_{w_+^\top - w_-^\top} w_+^\perp +
\Pi_+(w_+^\top, w_+^\top)\\
& + \Pi_-(w_-^\top, w_-^\top) + \nabla_{w_+} w_+ \cdot N_+ +
\nabla_{w_-} w_- \cdot N_- \} - \rho_+ w_+^\perp \nabla_{w_+}
p_{v,v}^+ - \rho_- w_-^\perp \nabla_{w_-} p_{v,v}^- \, dS
\end{split}\]
Using~\eqref{E:XYN} to compute $\nabla_{w_\pm} w_\pm \cdot N_\pm$,
we obtain
\[
|\int_{\RR^n \smallsetminus S_t} \rho \nabla p_{v,v} \nabla p_{w,w} dx|
\le C (|w|_{L^2 (S_t)}^2 + |w|_{L^2 (\RR^n \smallsetminus S_t)}^2).
\]
To compute $II_u(\bar v, \bar w)^2$, we use the decomposition
$p_{v,w}^\pm= \frac 1{\rho_\pm} \CH_\pm p_{v,w}^S -
\Delta_\pm^{-1}\text{tr} (DvDw)$ and the the Divergence Theorem to
obtain
\[
\int_{\RR^n \smallsetminus S_t} \rho |\nabla p_{v,w}|^2 dx = \int_{S_t}
p_{v,w}^S \CN p_{v,w}^S dS + \int_{\RR^n\smallsetminus S_t} \rho |\nabla
\Delta_\pm^{-1}\text{tr} (DvDw)|^2 dx.
\]
Therefore,
\[
|\int_{\RR^n \smallsetminus S_t} \rho |\nabla p_{v,w}|^2 dx - \int_{S_t}
p_{v,w}^S \CN p_{v,w}^S dS| \le C |w|_{L^2 (\RR^n \smallsetminus
S_t)}^2.
\]
Moreover, from expression~\eqref{E:pS} of $p_{w,v}^S$, we claim the
terms other than $\CN^{-1} \CD \cdot (w_\pm^\perp v_\pm^\top)$ are
of lower order. In fact, since
\[
|\nabla_{w_\pm} v_\pm - \nabla \Delta_\pm^{-1} \text{tr}
(DvDw)|_{L^2(\Omega_t^\pm)} \le C |w|_{L^2 (\Omega^\pm)}
\]
and it is divergence free, its normal component on $S_t$ is in
$H^{-\frac 12} (S_t)$ and we have
\[
|p_{v,w}^S + \CN^{-1} \CD \cdot (w_+^\perp (v_+^\top -
v_-^\top))|_{H^{\frac 12} (S_t)} \le C |w|_{L^2 (\RR^n \smallsetminus
S_t)}.
\]
Therefore, summarizing these estimates, we obtain
\begin{equation} \label{E:R}
|\bar \SR(u) (\bar v, \bar w)\bar v \cdot \bar w + \int_{S_t}
|\CN^{-\frac 12} \CD \cdot (w_+^\perp (v_+^\top - v_-^\top))|^2 dS |
\le C (|w|_{L^2 (S_t)}^2 + |w|_{L^2 (\RR^n \smallsetminus S_t)}^2).
\end{equation}
Using Lemma~\ref{L:SPrime}, we  define a self-adjoint
positive semi-~definite operator $\bar \SR_0(u) (\bar v)$ on $T_u
\Gamma$  which is the leading order part of $\bar \SR (u)(\bar v,
\cdot)\bar v$. In Eulerian coordinates, $\bar \SR_0(\bar v)$
takes the form
\begin{equation} \label{E:R0}\begin{cases}
\SR_0(v) (w) =\nabla f_+ \mathds{1}_{\Omega^+} + \nabla f_-\mathds{1}_{\Omega^-}, \\  f_\pm = \frac
1{\rho_+ \rho_-} \CH_\pm \CN^{-1} \CN_\mp \nabla_{v_+^\top -
v_-^\top} \CN^{-1} \CD \cdot (w_\pm^\perp (v_+^\top - v_-^\top)).\end{cases}
\end{equation}
Assuming smooth $S_t$ and $\bar v$, from~\eqref{E:R}, we can write
\[
\bar \SR (\bar v, \bar w) \bar w = \bar \SR_0 (\bar v) \bar w +
\text{ at most 1st order diff. operators}.
\]
Clearly $\bar \SR_0(\bar v)$ is a second order negative semi-definite
differential operator. Therefore, the linearized Euler's
equation~\eqref{E:criticalL} would be ill-posed if there had been no
surface tension, for $\bar \SR_0(\bar v)$ would become the leading
order term. {\it This is the Kelvin-Helmholtz instability of
vortex-sheets.}

Note that since $\CA$ is positive definite and is higher order than the self-adjoint $\SR_0(v)$, it is not difficult to see that the linearized  problem~\eqref{E:criticalL} is well-posed. For a priori estimates of the Euler's equation (E) with the boundary condition (BC), the positive semi-definiteness of the  leading order part
$\CA(u)$ of $\SD^2 S(u)$ suggests to consider the inner product
of~\eqref{E:critical} with $(\SD^2 S)^k u_t$ to obtain a priori estimates.

\section{Main Results} \label{S:estiE}

In this section, we will derive local energy estimate. We show that solutions of (E) with boundary
condition (BC) are locally bounded in
\begin{equation} \label{E:assump}
v(t, \cdot) \in H^{\frac32 k} (\RR^n \smallsetminus S_t) \quad \text{
and } \quad S_t \in H^{\frac32 k+1},
\end{equation}
where $k$ is an integer satisfying $\frac32 k> \frac n2 +1$
(equivalently $\frac32 k \ge \frac n2 + \frac 32$).

\noindent{\bf Definition of the energies and statements of the
theorems}. The conserved energy of  the Euler's equation is given by $E_0  = \int_{\RR^n \smallsetminus S_t} \frac 12 |v|^2 dx + S(S_t)$. Higher order
energies are based on the linearized Euler flow and thus involve the
differential operators $\CA$ defined in~\eqref{E:CA} and $\SD_t$
defined in~\eqref{E:barCD1}.

Let $\omega_v: \RR^n \to \RR^n$, often  written as $\omega$
for short, represent the curl a vector field $v$
defined on $\RR^n \smallsetminus S_t$, i.e.
\[
\omega (X) \cdot Y = \nabla_X v \cdot Y - \nabla_Y v \cdot X
\]
for any vector $X, Y \in \RR^n$. Viewing $\omega$ as a matrix, its
entries are $\omega_i^j = \omega(\frac \p{\p x^i}) \cdot \frac \p{\p
x^j} = \p_i v^j - \p_j v^i$.

\begin{defi} \label{E:energy}
Given  domains $\Omega^\pm$ with $\Omega_+$ compact and the
interface $S$ in $H^{\frac 32 k+1}$ and any vector field $v \in
H^{\frac 32 k}(\RR^n \smallsetminus S)$ with $v_+^\perp +v_-^\perp|_S=0$
and $\nabla \cdot v=0$, define the energy $E(S, v)$, often written
as $E$ for short,
\[
E = \int_{\RR^n \smallsetminus S} \frac 12 |\CA^{\frac k2} v|^2 + \frac
12 |\CA^{\frac k2 -\frac12} \nabla p_\kappa|^2\; dx +
|\omega|_{H^{\frac 32 k-1} (\RR^n \smallsetminus S)}^2,
\]
where $p_\kappa$ is the pressure due to the surface tension defined
in~\eqref{E:SPrime}.
\end{defi}
Since the free boundary is evolving, we  consider the following type of
$H^{\frac 32 k-\frac 12}$ neighborhoods of hypersurfaces  to maintain uniform constants in the energy  inequalities .
\begin{defi}\label{D:domainnbd}
Let $\Lambda= \Lambda (S, \frac 32 k-\frac 12, \delta)$ be the
collection of all hypersurfaces $\tilde S$  such that there
exists a diffeomorphism $F: S \to \tilde S\subset\RR^n$,  with $|F -
id_S|_{H^{\frac 32 k-\frac 12} (S)} < \delta$.
\end{defi}

Fix $0 < \delta \ll 1$ and let $\Lambda_0 \triangleq \Lambda (S_0,
\frac 32 k - \frac 12, \delta)$. From~\eqref{E:SPrime},
\eqref{E:CA}, and \eqref{E:CA1},
\begin{align}
&|p_\kappa|_{H^{s+\frac 12} (\RR^n \smallsetminus S)} \le C
|\kappa|_{H^s (S)}, \qquad &&s\in [\frac 12, \frac 32 k- \frac 12]
\label{E:pkappa}\\
&|\CA|_{L(H^s(\RR^n \smallsetminus S), H^{s-3}(\Omega))} \le C, \qquad
&&s \in [4-\frac 32 k, \frac 32 k-1]
\end{align}
where $C$ is uniform in $S \in \Lambda_0$. The next proposition
gives bounds on the velocity and mean curvature in terms of the
energy $E$.

\begin{proposition} \label{P:energy}
For $S \in \Lambda_0$ with $S \in H^{\frac 32 k+1}$, we have
\[
|\kappa|_{H^{\frac 32 k -1}(S)}^2 \le C_0 E, \qquad |v|_{H^{\frac 32
k} (\RR^n \smallsetminus S)}^2 \le C_0 (E+ E_0)^m
\]
for some integer $m>0$ depending only on $k$ and $n$ and some
constant $C_0>0$ depending only on the set $\Lambda_0$.
\end{proposition}

The proof of this proposition will be given below. Using this result
we will prove the following theorem on energy estimates.

\begin{theorem} \label{T:energy}
Fix $\delta>0$ sufficiently small. Then there exists $L>0$ such
that, if a solution of (E) and (BC) is given by $S_t$ with $S_t \in
H^{\frac 32 k+1}$ and $v(t, \cdot)\in C_t^0(H^{\frac 32 k}(\RR^n
\smallsetminus S_t))$, then there exists $t^*>0$, depending only on
$|v(0, \cdot)|_{H^{\frac 32 k} (\RR^n \smallsetminus S_t)}$, $L$, and
the set $\Lambda_0$, such that, for all $t \in [0, t^*]$,
\begin{equation}\begin{split}
&S_t \in \Lambda_0 \qquad \text{ and } \qquad |\kappa|_{H^{\frac 32
k-1} (S_t)} \le L, \\
&E(S_t, v(t, \cdot)) \le 2 E(S_0, v(0, \cdot)) + C_1 + \int_0^t P
(E_0, E(S_{t'}, v(t', \cdot))) \, dt' \label{E:energyE}
\end{split}
\end{equation}
where $P(\cdot)$ is a polynomial of positive coefficients determined
only by the set $\Lambda_0$ and $C_1$ is an constant determined only
by $|v(0, \cdot)|_{H^{\frac 32 k -\frac32}(\RR^n \smallsetminus S_0)}$,
and the set $\Lambda_0$.
\end{theorem}

Since the domain is evolving, the above continuity assumption of $v$
in $t$ means that there exist extensions of both $v_+$ and $v_-$ to
$[0, T] \times \RR^n$ which are continuous in $H^{\frac 32 k}(\RR^n)$.

In order to prove Proposition~\ref{P:energy} and
Theorem~\ref{T:energy}, we need the following lemmas.


\begin{lemma} \label{L:Pi} For any $S\in \Lambda_0$ with
$\kappa \in H^s (S)$, $s\in [\frac32 k-\frac 52, \frac32 k -1]$, we
have
\[
|\Pi|_{H^s (S)} + |N|_{H^{s+1} (S)} \le C(1+ |\kappa|_{H^s (S)}),
\]
for some $C>0$ uniform in $S \in \Lambda_0$.
\end{lemma}
\begin{proof}  The proof  is a straightforward  application of elliptic regularity applied to
\begin{equation} \label{E:DeltaPi}
- \Delta_S \Pi = - \CD^2 \kappa + (|\Pi|^2 I- \kappa \Pi) \Pi.
\end{equation}
and the fact that $|\Pi|_{H^{\frac n2 -1}(S)}\le C$ uniform in $S \in \Lambda_0$.
\end{proof}
\begin{coro} \label{C:nablaq} Suppose $S\in \Lambda_0$ with
$\kappa \in H^{\frac 32 k-\frac 32} (S)$, $g\in H^{\frac 32 k-1}
(\Omega^\pm)$, and $q = - \Delta_\pm^{-1}g$, then we have
\[
|\nabla_{N_{\CH_\pm}} q|_{H^{\frac 32 k} (\Omega^\pm)} \le C (1+
|\kappa|_{H^{\frac 32 k-\frac 32} (S)}) |g|_{H^{\frac 32 k-1}
(\Omega^\pm)}
\]
for some $C>0$ uniform in $S \in \Lambda_0$.
\end{coro}

The proof of this corollary follows Lemma~\ref{L:Pi} and the
identities
\begin{align*}
&\nabla_{N_\pm} \nabla_{N_{\CH_\pm}} q = {\CN(N_\pm)} \cdot \nabla
q+ D^2 q (N_\pm, N_\pm) =\CN(N_\pm)\cdot \nabla q -g  - \kappa_\pm
\nabla_{N_\pm} q \qquad &&\text{ on } S \\
&- \Delta \nabla_{N_{\CH_\pm}} q = \nabla_{N_{\CH_\pm}} g - 2 D^2 q
\cdot D N_{\CH_\pm} && \text{ in } \Omega^\pm.
\end{align*}

\begin{lemma} \label{L:CN} Suppose $S\in \Lambda_0$ with
$\kappa \in H^{\frac 32 k-\frac 32} (S)$,
\[
|(-\Delta_S)^{\frac12} - \CN_\pm|_{L(H^{s'} (S))} \le C (1+
|\kappa|_{H^{\frac 32 k-\frac 32} (S)}), \qquad s' \in [\frac
12-\frac 32 k, \frac 32 k-\frac 12].
\]
\end{lemma}

\begin{proof} From the identity
\[
(-\Delta_S - \CN_\pm^2) f = \kappa_\pm \CN_\pm(f) - 2\nabla_{N_\pm}
(- \Delta_\pm)^{-1} (D N_{\CH_\pm} \cdot D^2 f_{\CH_\pm}) - \CN_\pm
(N_\pm) \cdot (\CN_\pm (f) N_\pm + \nabla^\top f)
\]
for any smooth $f: S \to \RR$,  and lemma~\ref{L:Pi} we have
$$
|-\Delta_S - \CN_\pm^2|_{L(H^{s'} (S),
H^{s'-1} (S))}\le C
$$
Using commutators estimates $[-\Delta_S ,\CN_\pm]$ and factorization,  the lemma follows.  A detailed proof of the commutators and the factorization is given in section 6 of \cite{sz06}.
\end{proof}

\noindent {\bf Proof of Proposition~\ref{P:energy}.} The two terms
$|\CA^k v|_{L^2(S)}^2$ and $|\CA^{k -\frac12} \nabla
p_\kappa|_{L^2(S)}^2$ can be written explicitly using the
definition~\eqref{E:CA} of $\CA$
\begin{align}
&|\CA^{\frac k2} v|_{L^2(\RR^n \smallsetminus S)}^2 = \int_S v_+^\perp
( -\Delta_S \bar \CN )^{k-1} (-\Delta_S) v_+^\perp dS \label{E:energy1}\\
&|\CA^{\frac k2 -\frac12} \nabla p_\kappa|_{L^2(\RR^n \smallsetminus
S)}^2 = \int_S \kappa_+ \bar \CN ( -\Delta_S \bar \CN )^{k-1}
\kappa_+ dS, \label{E:energy2}
\end{align}
where
\begin{equation} \label{E:barN}
\bar \CN = (\frac 1{\rho_+} \CN_+) \CN^{-1} (\frac 1{\rho_-} \CN_-)
= ((\frac 1{\rho_+} \CN_+)^{-1} + (\frac 1{\rho_-}
\CN_-)^{-1})^{-1}.
\end{equation}
Clearly $\bar \CN$ is also self-adjoint and positive. The estimates
on $|\kappa|_{H^{\frac 32 k -1} (S)}$ and $|v_+^\perp|_{H^{\frac 32
k -\frac 12} (S)}$ follow immediately since from lemma~\ref{L:CN} $\CN_\pm$ behaves like
$(-\Delta_S)^{\frac 12}$.

To bound  $v$ in terms of $E_0$ and $E$ we note that $\Delta v$ is controlled  by $E$ from
\begin{equation} \label{E:deltav}
\Delta v^i = \p_j \omega_j^i.
\end{equation}

Therefore it is sufficient to control the boundary value $\nabla_{N_\pm} \, v_\pm$ by $E$. Moreover from the identity
$\nabla_{N_\pm} \, v_\pm =(D v_\pm)^*(N_\pm) + \omega_\pm (N_\pm)$ where $\omega_\pm$ is the restriction of $\omega$ on $S$,  it suffices to show that $E_0$ and $E$ control $\nu_\pm =
(D v_\pm)^*(N_\pm)$.


We first estimate $\nu_+^\top$ using the identity
\begin{equation} \label{E:Deltanu}
\Delta_S\, \nu_+^\top = \nabla^\top (\CD \cdot \nu_+^\top) +
\text{Ric}((D v_+)^*(N_+)^\top) +(\CD_{X_j}\, \omega_{\nu_+}^\top)
(X_j), \qquad \text{ at } x \in S
\end{equation}
where Ric is the Ricci curvature of $S$, $\{X_1, \ldots, X_{n-1}\}$
is any orthonormal frame of $T_x S$, and the tangential curl
$\omega_{\nu_+}^\top$ of $\nu_+^\top$ is defined as
\[
\omega_{\nu_+}^\top (x) (X) \cdot Y = \CD_X \, \nu_+^\top \cdot Y -
\CD_Y\, \nu_+^\top \cdot X = \nabla_X \, \nu_+^\top \cdot Y -
\nabla_Y \, \nu_+^\top \cdot X
\]
for any $X, Y \in T_x S$. From the definition of $\nu_\pm$ we have
\[
\omega_{\nu_+}^\top (X) \cdot Y = \Pi_+ (X) \cdot \nabla_Y \, v_+ -
\Pi_+(Y) \cdot \nabla_X \, v_+
\]
Therefore, by Sobolev inequalities, there exists $C>0$ uniform in $S
\in \Lambda_0$ so that
\[ \begin{split}
|\Delta_S \nu_+^\top|_{H^{\frac {3k -7}2} (S)} \le& |\CD \cdot
\nu_+^\top|_{H^{\frac {3k}2 - \frac 52} (S)} +
|\text{Ric}|_{H^{\frac {3k}2 -\frac 52} (S)} |N_+|_{H^{\frac {3k}2 -
\frac 32} (S)} |D v_+|_{H^{\frac {3k}2 -2} (S)} +
|\omega_{\nu_+}^\top|_{H^{\frac {3k}2 - \frac 52} (S)} \\
\le &|\CD \cdot \nu_+^\top|_{H^{\frac 32 k - \frac 52} (S)} + C
|v_+|_{H^{\frac32 k -\frac 18} (\Omega^+)}.
\end{split}\]
Here the norm $H^{\frac 32 k-\frac18}$ is chosen to illustrate that
the term is lower order. In fact any $H^{\frac 32 k -\alpha}$ with
$0 < \alpha < \frac 12$ works. To estimate the divergence term $\CD
\cdot \nu_+^\top$, one may compute
\[
\CD \cdot \nu_+^\top  = \Delta_S v_+^\perp -\CD \cdot (\Pi_+
(v_+^\top)),
\]
which along with Lemma~\ref{L:Pi} implies
\[\begin{split}
|\CD \cdot \nu_+^\top|_{H^{\frac {3k}2 - \frac 52} (S)} \le &C
|v_+^\perp|_{H^{\frac 32 k -\frac 12} (S)} + C |\Pi_+|_{H^{\frac 32
k -\frac 32} (S)} |v_+ - (v_+ \cdot N_+) N_+|_{H^{\frac 32 k -\frac
32} (S)}\\
\le &C|v_+^\perp|_{H^{\frac 32 k -\frac 12} (S)} + C (1 +
|\kappa_+|_{H^{\frac32 k -\frac 32} (S)}) |v_+|_{H^{\frac32 k -1}
(\Omega^+)}.
\end{split}\]
Therefore, from~\eqref{E:energy1} and~\eqref{E:energy2}, we obtain
\begin{equation} \label{E:Deltanu1}
|\Delta_S \nu_+^\top|_{H^{\frac 32 k-\frac 72} (S)} \le C E^{\frac
12} + C (1 + E^{\frac 13}) |v_+|_{H^{\frac32 k -\frac 18}
(\Omega^+)}.
\end{equation}

Finally, we only need to estimate $\nu_+^\perp =\nabla_{N_+} v_+
\cdot N_+$, much as in the way in~\cite{sz06}. Extending $\nu_+$
into $\Omega^+$ as $\nu_+ = (Dv_+)^* (N_{\CH_+})$, where $N_{\CH_+}$
is the harmonic extension of $N_+$ into $\Omega^+$, and comparing
the two ways of computing $\nabla \cdot \nu_+$ on $S$ using 1)
frames and $\omega$ and 2) divergence decomposition formula on $S$,
we obtain
\[
\nabla_{N_+}\, \nu_+ \cdot N_+ = (\nabla_{X_i} \omega)(X_i) \cdot
N_+ + D v_+ \cdot D N_{\CH_+} - \CD\cdot \nu_+^\top - \kappa_+
\nu_+^\perp.
\]
Moreover,
\[\begin{split}
\nabla_{N_+} \nu_+ \cdot N_+ =& \nabla_{N_+} (\nabla_{N_{\CH_+}} v_+
\cdot N_{\CH_+}) - N_+ \cdot D v_+(\CN_+ (N_{\CH_+})) \\
=& \CN_+ (\nu_+^\perp) + \nabla_{N_+} \Delta_+^{-1} \Delta
(\nabla_{N_{\CH_+}}\, v_+ \cdot N_{\CH_+}) - N_+ \cdot D v_+(\CN_+
(N_{\CH_+})).
\end{split}\]
Therefore, from the estimate on $\CD \cdot \nu_+^\top$, we obtain
\[
|\CN_+ \nu_+^\perp|_{H^{\frac {3k}2 - \frac 52} (S)} \le C E^{\frac
12} + C (1 + E^{\frac 13}) |v_+|_{H^{\frac32 k -\frac 18}
(\Omega^+)},
\]
which, along with~\eqref{E:deltav} and~\eqref{E:Deltanu1}, implies
\[
|v_+|_{H^{\frac 32 k} (\Omega^+)}^2 \le C E + C (1 + E^{\frac 23})
|v_+|_{H^{\frac32 k -\frac 18}\Omega^+)}^2.
\]
The estimate in Proposition~\ref{P:energy} follows immediately from
Sobolev inequalities. \hfill $\square$\\

\noindent {\bf Proof of Theorem~\ref{T:energy}.} To prove
Theorem~\ref{T:energy}, in addition to Proposition~\ref{P:energy},
we need the following: a) the estimates on the Lagrangian
coordinates map and consequently $\kappa\in{H^{\frac 32 k-\frac 52}
(S)}$, b) estimates on $\omega = D v - (D v)^*$, and  c) commutators
involving $\D_t$.  In the following all constant $C>0$ will depend only on the set $\Lambda_0$.\\


\noindent {\it Estimate of the Lagrangian coordinate map $u(t, y)$.}
We will only work on the domain $\Omega_t^t$. From our assumption on
$v$, the ODE $u_t (t, y)= v(t, u(t, y))$ solving $u$ is well-posed.
Since $u(t, \cdot): \Omega_0^+ \to \Omega_t^+$ is volume preserving
and $\frac 32k>\frac n2 +1$, it is easy to derive, Therefore,
\begin{equation}\label{E:LagE1}
|u(t, \cdot)-I|_{H^{\frac 32 k} (\Omega_0^+)} \le C \int_0^t |v(t',
\cdot)|_{H^{\frac 32 k} (\Omega_t^+)} |u(t', \cdot)|_{H^{\frac 32 k}
(\Omega_0^+)}^{\frac 32 k} \;dt',
\end{equation}
where $C>0$ depends only on $n$ and $k$. Let $\mu>0$ be a positive
large number specified later,
\begin{equation} \label{E:t0}
t_0 = \sup \{t \mid |v(t',\cdot)|_{H^{\frac 32 k} (\RR^n \smallsetminus
S_t)} < \mu, \; \forall t' \in [0, t]\},
\end{equation}
We have $t_0>0$ due to the continuity of $v(t, \cdot)$ in $H^{\frac
32 k} (\Omega_t)$. From ODE estimates, there exists $t_1>0$ and
$C_2>0$ which depend only on $\mu$ such that, for all $0\le t \le
\min\{t_0, t_1\}$,
\begin{equation}\label{E:LagE}
|u(t, \cdot)-I|_{H^{\frac 32 k} (\Omega_0^+)} \le C_2 t.
\end{equation}
It implies the mean curvature estimate, for all $0\le t \le
\min\{t_0, t_1\}$,
\begin{equation}
|\kappa(t, \cdot)|_{H^{\frac 32 k-\frac 52} (S_t)} \le |\kappa(0,
\cdot)|_{H^{\frac 32 k-\frac 52} (S_0)} + C_3 t.
\end{equation}
Here it is easy to see from local coordinates that $C_3$ is
determined only by $\mu$ and the set $\Lambda_0$. Therefore, there
exists $t_2>0$ determined only by $\mu$ and the set $\Lambda_0$ such
that $S_t \in \Lambda_0$ for $0\le t \le \min\{t_0, t_2\}$.\\

\noindent {\it Evolution of the curl $\omega = D v - (Dv)^*$.} It is
easy to compute
\[
\D_t \omega = D \D_t v - (D \D_t v)^* + ((Dv)^*)^2 - (Dv)^2 =
((Dv)^*)^2 - (Dv)^2 = - (Dv)^* \omega - \omega Dv.
\]
Since $\frac 32k > \frac n2 +1$, we have
\begin{equation} \label{E:dtomega}
\frac d{dt} \int_{\RR^n \smallsetminus S_t} |D^{\frac 32 k -1}
\omega|_{H^{\frac 32 k-1} (\RR^n \smallsetminus S_t)}^2 dx \le C
|v|_{H^{\frac 32 k} (\RR^n \smallsetminus S_t)} |\omega|_{H^{\frac 32
k-1} (\RR^n \smallsetminus S_t)}^2.
\end{equation}

\noindent {\it The commutator involving $\D_t$.} First we list a few
commutators calculated in~\cite{sz06}:
\begin{align}
&\D_t \nabla g = \nabla \D_t g - (D v)^* (\nabla g),
\label{E:dtnabla}\\
&[\D_t, \Delta_\pm^{-1}] g = \Delta_\pm^{-1} ( 2D v \cdot D^2
\Delta_\pm^{-1} g + \Delta v \cdot \nabla \Delta_\pm^{-1} g),
\label{E:dtdelta-1}
\end{align}
for any smooth function $g$. It also has been proved in~\cite{sz06}
that, for any function $f$ defined on $S_t$,
\begin{align}
& |[\D_t, \Delta_{S_t}]|_{L(H^{s_1} (S_t), H^{s_1-2} (S_t))} \le C
|v|_{H^{\frac 32 k} (\RR^n \smallsetminus S_t)} \qquad &&s_1 \in
(\frac 72-\frac 32 k, \frac 32 k -\frac12], \label{E:dtdeltat1}\\
& |[\D_t, \CN_\pm]|_{L(H^{s_2} (S_t), H^{s_2-1}(S_t))} \le C
|v|_{H^{\frac 32 k} (\RR^n \smallsetminus S_t)}  &&s_2 \in [\frac 12,
\frac 32 k -\frac12]. \label{E:dtcnW1}
\end{align}
It also implies that
\begin{align}
&|[\D_t, \CN^{-1}]|_{L(H^s (S_t), H^{s-1}(S_t))} \le C |v|_{H^{\frac
32 k} (\RR^n \smallsetminus S_t)}  \qquad &&s \in [-\frac 12, \frac 32 k
-\frac32] \label{E:dtcnW2}\\
&|[\D_t, \bar \CN]|_{L(H^s (S_t), H^{s-1}(S_t))} \le C |v|_{H^{\frac
32 k} (\RR^n \smallsetminus S_t)}  &&s \in [\frac 12, \frac 32 k
-\frac12] \label{E:dtbarcn}
\end{align}

\noindent {\it Evolution of $E$.} For the rest of this section, let
$Q = Q(|v|_{H^{\frac 32 k} (\RR^n \smallsetminus S_t)},
|\kappa|_{H^{\frac 32 k-1} (S_t)})$ denote a generic positive
polynomial in $|v|_{H^{\frac 32 k} (\RR^n \smallsetminus S_t)}$,
and $|\kappa|_{H^{\frac 32 k-1} (S_t)}$ with coefficients depending
only on the set $\Lambda_0$. One may notice here, by

Lemma~\ref{L:Pi}, $|\Pi|_{H^{\frac 32 k -1} (S_t)}$, and
$|v_\pm^\perp|_{H^{\frac 32 k -\frac 12} (S_t)}$, and
$|v_\pm^\top|_{H^{\frac 32 k -\frac 12} (S_t)}$ can also be included
in $Q$. Recall that
\[
\D_{t_\pm} dS = (\CD \cdot v_\pm^\top + \kappa_\pm v_\pm^\perp) dS.
\]
Since
\[
\nabla \cdot v_\pm|_{S_t} = \CD \cdot v_\pm^\top + \kappa_\pm
v_\pm^\perp + \nabla_{N_\pm} v_\pm \cdot N_\pm=0
\]
then
\[
|\kappa_\pm v_\pm^\perp + \CD \cdot v_\pm^\top|_{H^{\frac 32 k-\frac
32} (S_t)} = |\nabla_{N_\pm} v_\pm \cdot N_\pm|_{H^{\frac 32 k-\frac
32} (S_t)} \le C |v|_{H^{\frac 32 k} (\RR^n \smallsetminus S_t)}
\]
and thus $\D_t dS$ would not complicate the estimates since $\frac
32 k\ge \frac n2 +\frac 32$.

Before we embark on calculating the energy inequality it is helpful to recall two facts. First,   we only need to keep track of terms which can not be bounded by $|v_\pm|_{H^{\frac 32 k}(\Omega^\pm_t)}$ and $|\kappa|_{H^{\frac 32 k -1}(S_t)}$, and second, $\CN_\pm$ are selfadjoint operators of order one.

\noindent {\bf I:} $\qquad \lf|\frac 12 \frac d{dt} |\CA^{\frac k2
-\frac12} \nabla p_\kappa|_{L^2(\RR^n \smallsetminus S_t)}^2 -
\int_{S_t} \kappa_+ \bar \CN ( -\Delta_{S_t} \bar \CN )^{k-1}
(-\Delta_{S_t}) v_+^\perp
dS\rt| \le Q.$ \\
From~\eqref{E:energy2}, \eqref{E:dtdeltat1}, and~\eqref{E:dtbarcn},
it is clear
\[
\lf|\frac 12 \frac d{dt} |\CA^{\frac k2 -\frac12} \nabla
p_\kappa|_{L^2(\RR^n \smallsetminus S_t)}^2 - \int_{S_t} \kappa_+ \bar
\CN ( -\Delta_{S_t} \bar \CN )^{k-1} \D_{t_+} \kappa_+ dS\rt| \le Q.
\]
From the expression for $\D_{t_\pm} \kappa_\pm $ given in \cite{sz06}
\begin{equation}
\D_{t_\pm} \kappa_\pm = -\Delta_{S_t} v_\pm^\perp - v_\pm^\perp
|\Pi|^2 + (\CD \cdot \Pi_\pm)(v_\pm^\top)= -\Delta_{S_t} v_\pm^\perp
- v_\pm^\perp |\Pi|^2 + \nabla_{v_\pm^\top} \kappa_\pm \label{E:dtk}
\end{equation}
and Lemma~\ref{L:Pi} we only need
\begin{equation} \label{E:dtk1}
|\int_{S_t} \kappa_+ \bar \CN ( -\Delta_{S_t} \bar \CN )^{k-1}
\nabla_{v_+^\top} \kappa_+ dS | \le Q.
\end{equation}
%
to derive estimate {\bf I}.  By considering a flow $\phi(\tau, \cdot)$ on $\Omega_t^+$ generated
by $ \CH_+ v_+^\top$, the above commutator estimates applied to $\D_\tau$ allow us  to pull $\nabla_{v_+^\top}$ to the front with lower order errors
\[
\lf|\int_{S_t} \kappa_+ \bar \CN ( -\Delta_{S_t} \bar \CN )^{k-1}
\nabla_{v_+^\top} \kappa_+ dS - \frac 12 \int_{S_t}
\nabla_{v_+^\top} \lf(\kappa_+ \bar \CN ( -\Delta_{S_t} \bar \CN
)^{k-1} \kappa_+ \rt) dS \rt| \le Q.
\]
Therefore, inequality~\eqref{E:dtk1} follows from the Divergence
Theorem and {\bf I} follows consequently.\\

\noindent {\bf II}: $\qquad \lf|\frac d{dt} (\frac 12 |\CA^{\frac
k2} v|_{L^2(\RR^n \smallsetminus S_t)}^2 - E_{ex}) + \int_{S_t}
v_+^\perp ( -\Delta_{S_t} \bar \CN )^k \kappa_+ \,dS \rt| \le Q$ \\
where the extra term
\[ \begin{split}
E_{ex} =  \frac {\rho_+}{2(\rho_+ + \rho_-)} \int_{S_t}
\nabla_{v_+^\top} \kappa_+ \cdot &\bar \CN (-\Delta_{S_t} \bar \CN
)^{k-2} \nabla_{v_+^\top} \kappa_+ \, dS \\
&-  \frac {\rho_-}{2(\rho_+ + \rho_-)} \int_{S_t} \nabla_{v_-^\top}
\kappa_+\cdot \bar \CN ( -\Delta_{S_t} \bar \CN )^{k-2}
\nabla_{v_-^\top} \kappa_+ \, dS \}.
\end{split} \]
From~\eqref{E:energy1}, \eqref{E:dtdeltat1}, and~\eqref{E:dtbarcn},
it is clear
\begin{equation} \label{E:v1}
\lf|\frac 12 \frac d{dt} |\CA^{\frac k2} v|_{L^2(\RR^n \smallsetminus
S_t)}^2 - \int_{S_t} v_+^\perp ( -\Delta_{S_t} \bar \CN )^{k-1}
(-\Delta_{S_t}) \D_{t_+} v_+^\perp dS \rt| \le Q.
\end{equation}
Using~\eqref{E:euler}, \eqref{E:dtn}, \eqref{E:SPrime},
and~\eqref{E:II}, we have
\[ \begin{split}
\D_{t_+} v_+^\perp =& (\D_{t_+} v_+) \cdot N_+ + v_+\cdot \D_{t_+}
N_+ = - \nabla_{N_+} p_{v,v}^+ - \nabla_{N_+} p_\kappa^+ -
\nabla_{v_+^\top} v_+ \cdot N_+\\
= & -\frac 1{\rho_+} \CN_+ p_{v,v}^S - \nabla_{N_+} \Delta_+^{-1}
\text{tr}(Dv)^2 - \bar \CN \kappa_+ - \nabla_{v_+^\top} v_+^\perp +
\Pi_+(v_+^\top, v_+^\top).
\end{split}\]
From the corollary of Lemma~\ref{L:Pi}, we have
\begin{equation} \label{E:pvv}
|\nabla_{N_\pm} \Delta_\pm^{-1} \text{tr}(Dv)^2|_{H^{\frac 32 k -
\frac 12} (S_t)} \le Q.
\end{equation}
Therefore,
\begin{equation} \label{E:dtvperp1} \begin{split}
|\frac 12 \frac d{dt} |\CA^{\frac k2} v|_{L^2(\RR^n \smallsetminus
S_t)}^2 - \int_{S_t}  v_+^\perp ( -\Delta_{S_t} \bar \CN )^{k-1}
(-\Delta_{S_t}) (& -\frac 1{\rho_+} \CN_+ p_{v,v}^S - \bar \CN
\kappa_+ \\
&- \nabla_{v_+^\top} v_+^\perp + \Pi_+(v_+^\top, v_+^\top) ) dS |
\le Q.
\end{split}\end{equation}
From~\eqref{E:II}, we can write the first term in the second
integral as
\[ \begin{split}
-\frac 1{\rho_+} \CN_+ p_{v,v}^S = \frac 1{\rho_+} \CN_+ \bar
\CN^{-1} \{ 2 \nabla_{v_+^\top - v_-^\top} v_+^\perp -
&\Pi_+(v_+^\top, v_+^\top) - \Pi_-(v_-^\top, v_-^\top) \\
&- \nabla_{N_+} \Delta_+^{-1} \text{tr}(Dv)^2 - \nabla_{N_-}
\Delta_-^{-1} \text{tr}(Dv)^2\}
\end{split} \]
Since $|\bar \CN^{-1}\CN_+|_{L(H^{\frac 32 k -\frac 12} (S_t))} \le
C$, we have that its dual $\CN_+\bar \CN^{-1}$ satisfies the same
estimate. Therefore, using~\eqref{E:pvv}, we obtain
\[
|-\frac 1{\rho_+} \CN_+ p_{v,v}^S - \frac 1{\rho_+} \CN_+ \bar
\CN^{-1} (2 \nabla_{v_+^\top - v_-^\top} v_+^\perp - \Pi_+(v_+^\top,
v_+^\top) - \Pi_-(v_-^\top, v_-^\top))|_{H^{\frac 32 k -\frac 12}
(S_t)} \le Q.
\]
From Lemma~\ref{L:CN}, we have
\[
|\CN_+ \bar \CN^{-1} - \frac {\rho_+\rho_-} {\rho_+ +
\rho_-}|_{L(H^{\frac 32 k -\frac 32} (S_t), H^{\frac 32 k -\frac 12}
(S_t))} \le Q,
\]
which, along with Lemma~\ref{L:Pi}, implies
\[
|-\frac 1{\rho_+} \CN_+ p_{v,v}^S - \frac {\rho_-} {\rho_+ +
\rho_-}(2 \nabla_{v_+^\top - v_-^\top} v_+^\perp - \Pi_+(v_+^\top,
v_+^\top) - \Pi_-(v_-^\top, v_-^\top))|_{H^{\frac 32 k -\frac 12}
(S_t)} \le Q.
\]
Substituting this estimate into~\eqref{E:dtvperp1}, we have
\[ \begin{split}
|\frac 12 \frac d{dt} |\CA^{\frac k2} v|_{L^2(\RR^n \smallsetminus
S_t)}^2 - \int_{S_t} & v_+^\perp ( -\Delta_{S_t} \bar \CN )^{k-1}
(-\Delta_{S_t}) \{ \frac {\rho_+}{\rho_+ + \rho_-} \Pi_+(v_+^\top,
v_+^\top) - \frac {\rho_-}{\rho_+ + \rho_-}
\Pi_-(v_-^\top, v_-^\top) \\
& - \bar \CN \kappa_+ + \nabla v_+^\perp \cdot (\frac {\rho_- -
\rho_+} {\rho_+ + \rho_-} v_+^\top - \frac {2\rho_-} {\rho_+ +
\rho_-} v_-^\top )\} dS | \le Q.
\end{split} \]
Much as the proof of inequality~\eqref{E:dtk1}, by commuting $\nabla_{v^\top_\pm}$ we  obtain
\[
|\int_{S_t} v_+^\perp ( -\Delta_{S_t} \bar \CN )^{k-1}
(-\Delta_{S_t}) \nabla v_+^\perp \cdot (\frac {\rho_- - \rho_+}
{\rho_+ + \rho_-} v_+^\top - \frac {2\rho_-} {\rho_+ + \rho_-}
v_-^\top ) dS| \le Q.
\]
Therefore, we have
\begin{equation} \label{E:dtvperp2}\begin{split}
|\frac 12 \frac d{dt} |\CA^{\frac k2} v|_{L^2(\RR^n \smallsetminus
S_t)}^2 - \int_{S_t} v_+^\perp ( -\Delta_{S_t} \bar \CN )^{k-1} &
(-\Delta_{S_t}) \{ \frac {\rho_+}{\rho_+ + \rho_-} \Pi_+(v_+^\top,
v_+^\top) \\
&- \frac {\rho_-}{\rho_+ + \rho_-} \Pi_-(v_-^\top, v_-^\top)  - \bar
\CN \kappa_+ \} dS | \le Q.
\end{split} \end{equation}
In order to estimate the terms with $\Pi$, we first use
Lemma~\ref{L:Pi} and identity~\eqref{E:DeltaPi} to obtain
\[
|\Delta_{S_t} (\Pi_\pm (v_\pm^\top, v_\pm^\top)) - \CD^2 \kappa_\pm
(v_\pm^\top, v_\pm^\top)|_{H^{\frac 32 k - \frac 52} (S_t)} \le Q.
\]
Since
\[
\CD^2 \kappa_\pm (v_+^\top, v_+^\top) = \nabla_{v_\pm^\top}
\nabla_{v_\pm^\top} \kappa_\pm - \CD_{v_\pm^\top} v_\pm^\top \cdot
\nabla \kappa_\pm
\]
we have
\[\begin{split}
|\frac 12 \frac d{dt} |\CA^{\frac k2} v|_{L^2(\RR^n \smallsetminus
S_t)}^2 - \int_{S_t} v_+^\perp ( -\Delta_{S_t} \bar \CN )^{k-1} \{
&\Delta_{S_t} \bar \CN \kappa_+ - \frac {\rho_+}{\rho_+ + \rho_-}
\nabla_{v_+^\top} \nabla_{v_+^\top} \kappa_+ \\
& + \frac {\rho_-}{\rho_+ + \rho_-} \nabla_{v_-^\top}
\nabla_{v_-^\top} \kappa_- \} dS| \le Q.
\end{split} \]
As in the proof of\eqref{E:dtk1}, we apply the commutator
estimates~\eqref{E:dtdeltat1} and~\eqref{E:dtcnW1} to move one of
$\nabla_{v_\pm^\top}$ and obtain
\[ \begin{split}
|\frac 12 &\frac d{dt} |\CA^{\frac k2} v|_{L^2(\RR^n \smallsetminus
S_t)}^2 - \frac {\rho_+}{\rho_+ + \rho_-} \int_{S_t}
\nabla_{v_+^\top} (-\Delta_{S_t} v_+^\perp) \cdot \bar \CN (
-\Delta_{S_t} \bar \CN )^{k-2} \nabla_{v_+^\top} \kappa_+ \, dS \\
&- \frac {\rho_-}{\rho_+ + \rho_-} \int_{S_t} \nabla_{v_-^\top} (
-\Delta_{S_t} v_+^\perp)\cdot \bar \CN ( -\Delta_{S_t} \bar \CN
)^{k-2} \nabla_{v_-^\top} \kappa_+ \, dS + \int_{S_t} v_+^\perp (
-\Delta_{S_t} \bar \CN )^k \kappa_+ \,dS | \le Q.
\end{split} \]
From identity~\eqref{E:dtk},
\[
|-\Delta_{S_t} v_+^\perp - \D_{t_+} \kappa_+|_{H^{\frac 32 k - 2}
(S_t)} \le Q,
\]
which implies
\[ \begin{split}
|\frac 12 &\frac d{dt} |\CA^{\frac k2} v|_{L^2(\RR^n \smallsetminus
S_t)}^2 - \frac {\rho_+}{\rho_+ + \rho_-} \int_{S_t}
\nabla_{v_+^\top} \D_{t_+} \kappa_+ \cdot \bar \CN
(-\Delta_{S_t} \bar \CN )^{k-2} \nabla_{v_+^\top} \kappa_+ \, dS \\
&- \frac {\rho_-}{\rho_+ + \rho_-} \int_{S_t} \nabla_{v_-^\top}
\D_{t_+} \kappa_+\cdot \bar \CN ( -\Delta_{S_t} \bar \CN )^{k-2}
\nabla_{v_-^\top} \kappa_+ \, dS + \int_{S_t} v_+^\perp (
-\Delta_{S_t} \bar \CN )^k \kappa_+ \,dS | \le Q.
\end{split} \]
Using~\eqref{E:dtdeltat1} and~\eqref{E:dtcnW1} one more time, we
obtain inequality {\bf II}. \\

\noindent {\it Proof of Theorem~\ref{T:energy}.} Adding
inequalities~\eqref{E:dtomega}, {\bf I}, and {\bf II}, we have
\begin{equation} \label{E:energy3}
E(t) -E(0) + E_{ex}(t) - E_{ex} (0) \le \int_0^t Q(|v|_{H^{3k}
(\RR^n \smallsetminus S_{t'})}, |\kappa|_{H^{3k-1} ( S_{t'})})\; dt',
\end{equation}
where $Q$ is a polynomial with positive coefficients that depends
only on $\Lambda_0$. This inequality holds on $[0,  \min\{t_0,
t_2\}]$ where $t_0$ is defined in~\eqref{E:t0} and $t_2$ is
determined only by $|v(0, \cdot)|_{H^{3k}(\RR^n \smallsetminus S_t)}$
and the set $\Lambda_0$, the neighborhood of $S_0$ in $H^{\frac 32 k
-\frac 12}$. Clearly
\[
|E_{ex}| \le C |v|^2_{H^{\frac 32 k-\frac58} (\RR^n \smallsetminus S_t)}
|\kappa|_{H^{\frac 32 k - \frac 32}(S_t)}^2,
\]
where $C>0$ depends only on the set $\Lambda_0$. Interpolating $v$
between $H^{3k}(\RR^n \smallsetminus S_t)$ and $H^{3k-\frac 32}(\RR^n
\smallsetminus S_t)$ and $\kappa$ between $H^{3k-1}(S_t)$ and
$H^{3k-\frac 52} (S_t)$, we obtain from Proposition~\ref{P:energy},
\[
|E_{ex}| \le \frac 12 E + C_1 (1+ |v|_{H^{3k-\frac32} (\RR^n
\smallsetminus S_t)}^m)
\]
for some integer $m>0$ where the constant $C_1$, which include
$|\kappa|_{H^{3k -\frac 52}(S_t)}$, is determined only by $E_0$ and
the set $\Lambda_0$. From the Euler's equation~\eqref{E:euler},
\eqref{E:pkappa}, \eqref{E:II}, Lemma~\ref{L:Pi}, and its corollary,
\[
|\D_t v|_{H^{\frac 32 k -\frac 32} (\RR^n \smallsetminus S_t)} = |\nabla
p_{v, v} +\nabla p_\kappa|_{H^{\frac 32 k -\frac 32} (\RR^n
\smallsetminus S_t)} \le Q.
\]
We can use the Lagrangian coordinate map $u(t, \cdot)$ to estimate
\[
|v(t, \cdot)|_{H^{3k -\frac 32} (\RR^n \smallsetminus S_t)} - |v(0,
\cdot)|_{H^{3k -\frac 32} (\RR^n \smallsetminus S_0)}
\]
Through a similar procedure of the derivation of~\eqref{E:LagE},
there exists $t_3>0$, depending only on $|v(0,\cdot)|_{H^{3k}(\RR^n
\smallsetminus S_t)}$ and the set $\Lambda_0$ so that for $0 \le t \le
\min\{t_0, t_3\}$,
\[
\lf| |v(t, \cdot)|_{H^{3k -\frac 32} (\RR^n \smallsetminus S_t)}^m -
|v(0, \cdot)|_{H^{3k -\frac 32} (\RR^n \smallsetminus S_0)}^m \rt| \le
\int_0^t Q \, dt'
\]
for some polynomial  $Q$ with positive coefficients. Therefore,
\[
E_{ex} \le \frac 12 E + C_1 (1+ |v(0, \cdot)|_{H^{3k -\frac 32}
(\RR^n \smallsetminus S_0)}^m) + \int_0^t Q  \, dt' \le  \frac 12 E +
C_1 + \int_0^t Q \, dt',
\]
where $C_1$ is determined only by $|v(0, \cdot)|_{H^{\frac32 k
-\frac32}(\RR^n \smallsetminus S_t)}$ and the set $\Lambda_0$. Thus
\[
E(S_t, v(t, \cdot)) \le 2 E(S_0, v(0, \cdot)) + C_1 + \int_0^t Q \,
dt'.
\]
By inserting the above inequality into~\eqref{E:energy3} and using
proposition~\ref{P:energy}, we obtain~\eqref{E:energyE}. By choosing
$\mu$ large enough compared to the initial data,
Theorem~\ref{T:energy} follows. \hfill $\square$
\\

\noindent{\bf  Notation}

\noindent
 \\  $A^*$: the adjoint operator of an operator.
 \\  $D$ and $\p$: differentiation with respect to spatial variables.
 \\  $\nabla f$: the gradient vector of a scalar function $f$.
 \\  $\nabla_X$:  the directional directive in the direction $X$.
 \\  $\perp$ and $\top$: the normal and the tangential components of the relevant quantities.
 \\  $\D_t = \p_t + v^i\p_{x^i}$ the material derivative  along the particle path.
 \\$S_t= \p\Omega_t$ the boundary of a smooth domain evolving in time.
\\  $N(t,x)$: the outward unit normal vector of $S_t$ at $x \in S_t$.
 \\  $\Pi$: the second fundamental form of $S_t$, $\Pi(t, x)(w) = \nabla_w N \in T_x S_t$.
\\ $\Pi(X, Y)=\Pi(X) \cdot Y$.
\\  $\kappa$: the mean curvature of $S_t$, i.e. $\kappa = \text{tr} \Pi$.
\\$f_\CH=\CH(f)$: the harmonic extension of $f$ on $\Omega_t$.
\\$\CN (f) = \nabla_N \CH(f) : S \to \RR$:  the Dirichlet-Neumann operator.
\\$\bar X = X\circ u^{-1}$ the Lagrangian coordinates description of $X$.
\\  $\CD$: the covariant differentiation on $ S_t \subset \RR^n$.
\\$\CD_w = \nabla_w^\top$,  for any $ x \in S_t \quad w\in T_x S_t$.
 \\  $\R(X, Y)$, $\; X, Y\in T_x  S_t$: the curvature tensor of $ S_t$.
 \\  $\Delta_\CM \triangleq \text{tr} \CD^2$: the Beltrami-Lapalace operator on a Riemannian manifold $\CM$.
 \\$\Delta^{-1}$: the inverse Laplacian with zero Dirichlet data.
 \\$\Gamma = \{\Phi_\pm \, : \Omega^\pm_t \to \RR^n \text{ ; volume preserving homeomorphism, such that } \, \Phi_+(\p\Omega^+)  = \Phi_-(\p\Omega^-)  \}$
 \\  $\bar \SD$: the covariant derivative on $\Gamma$,
 \\  $ \SD$: represent $\bar\SD$ in Eulerian coordinates.
 \\  $\bar \SR$: the curvature operator on $\Gamma$.
 \\  $ \SR$: represent $\bar\SR$ in Eulerian coordinates.
 \\  $\text{II}$: the second fundamental form of $\Gamma \subset L^2$
\\$\text{II}_u(w_1, w_2) = \nabla^\perp_{w_1}w_2$,  for any $u \in \Gamma, \quad w_1, w_2 \in T_u \Gamma$
\\ $p_{v,w} =-\Delta^{-1} \text{tr} (DvDw)$.

\end{document}